\documentclass{amsart}
\usepackage[applemac]{inputenc}
\usepackage{amsfonts,amsthm,latexsym,amsmath,amssymb,amscd,epsfig}
\usepackage{graphics,graphicx}

\begin{document}

\bigskip

\begin{center}
{\bf\Large Subharmonic Configurations and \\[1ex]
Algebraic Cauchy Transforms of Probability Measures}
\end{center}

\bigskip

\bigskip

\begin{center}
{\sc Jan-Erik Bj\"ork, Julius Borcea, Rikard B\o gvad}
\end{center}

\bigskip

\begin{center}
Department of Mathematics, Stockholm University, SE-106 91 Stockholm, Sweden\\
{\em E-mail addresses}: 
\texttt{jeb@math.su.se, julius@math.su.se, rikard@math.su.se}
\end{center}

\bigskip

\bigskip

\centerline{\bf Abstract}
\begin{quote}
We study subharmonic functions whose 
Laplacian is supported on a null set $K\subset \mathbf{C}$ and in 
connected components 
of $\mathbf{C}\setminus K$ admit harmonic extensions to larger sets. We prove 
that if such a function has a piecewise holomorphic derivative then it is 
locally piecewise harmonic and in generic cases it coincides locally with 
the maximum of finitely many harmonic functions. Moreover, we describe $K$ 
when the holomorphic derivative 
satisfies a global algebraic equation. The proofs follow 
classical patterns and our methods may also 
be of independent interest.
\end{quote}

\bigskip

\bigskip

\begin{center}
{\bf\large Introduction}
\end{center}

\bigskip

\noindent
Let $\Omega$ be an open connected
subset of the complex plane $\mathbf{C}$.
Denote by $\text{SH}_0(\Omega)$ the class of subharmonic functions 
$V$ in $\Omega$ for which the support of the Laplacian $\Delta(V)$
has Lebesgue measure $0$, where $\Delta(V)$ in the
sense of distributions  is a non-negative Riesz measure supported by
the null set
$\text{supp}(\Delta(V))$.
As explained in H\"ormander \cite{Horm}, every 
$V\in \text{SH}_0(\Omega)$
is identified with an element in
$L^1_{\text{loc}}(\Omega)$
and can always be taken as
an upper semi-continuous function. Moreover, the distribution derivatives
$\partial V/\partial x$ and $\partial V/\partial y$
belong to $L^1_{\text{loc}}(\Omega)$.
In particular, the distribution derivative
\[
\partial V/\partial z
=\frac{1}{2}\bigl(\partial V/\partial x-\partial V/\partial y \bigr)
\]
is a holomorphic function in
$\Omega\setminus\text{supp}(\Delta(V))$
which as a distribution is an element of
$L^1_{\text{loc}}(\Omega)$.
Therefore, if
the holomorphic function
$\partial V/\partial z$ defined in 
$\Omega\setminus\text{supp}(\Delta(V))$
extends to a holomorphic function
$g$ defined in the whole set
$\Omega$, then the distribution $\partial/\partial \bar z(\partial V/\partial z)= \Delta V/4=0$, i.e., $V$ is harmonic in $\Omega$.

\medskip

\noindent
The aforementioned facts, already known to
F.~Riesz who laid the foundations of subharmonic functions in his
famous article \cite{Riesz} from 1926, have led to
the problems studied in the present paper.
We call $V\in\text{SH}_0(\Omega)$ \emph{piecewise harmonic}
if there exists a finite set of harmonic functions $H_1,\ldots,H_k$
in $\Omega$ such that for every connected component $U$
of
$\Omega\setminus\text{supp}(\Delta(V))$ one has 
$V=H_j$ in $U$ for some $1\leq j\leq k$.
In this case we refer to  $V$ as a {\em subharmonic configuration} of the 
$k$-tuple $H_1,\ldots,H_k$. When $V$ is such a subharmonic configuration
one easily shows the inclusion
\[ 
\text{supp}(\Delta(V))\subset
\,\bigcup_{i\neq j}\,\{H_i=H_j\}.
\]

\noindent
In general 
a $k$-tuple of harmonic functions
$H_1,\ldots,H_k$  gives rise to \emph{several}
subharmonic configurations, see \S 2.11. 
An obvious subharmonic configuration  is  the 
maximum function
$V^*=\text{max}(H_1,\ldots,H_k)$. 
In Theorem 1.4 we show that
$V^*$ is locally the \emph{unique} subharmonic configuration 
of $H_1,\ldots,H_k$ in a neighborhood of a point
$p\in\Omega$
when
the $k$-tuple of gradient vectors
$\nabla(H_1)(p),\ldots,\nabla(H_k)(p)$ are extreme points of
their convex hull.
An essential role in proving this as well as our
other results 
is played by the \emph{Key Lemma 1.1} in Section 1.
\medskip

\noindent
The next issue in this article is to study
functions $V$ in $\text{SH}_0(\Omega)$
for which the analytic function in
$\Omega\setminus\text{supp}(\Delta(V))$
defined by
$\partial V/\partial z$
is \emph{piecewise holomorphic}.
This  means that there exists a finite set of 
holomorphic functions $g_1,\ldots,g_k$ in $\Omega$ and for 
every connected subset $U$ of 
$\Omega\setminus\text{supp}(\Delta(V))$ some
$1\leq j\leq k$ such that
$\partial V/\partial z=g_j$ in $U$.
Let us remark that if
$V$ is piecewise harmonic with respect to 
$H_1,\ldots,H_k$ then
$\partial V/\partial z$ is piecewise holomorphic
with respect to the $k$-tuple $\{\partial H_i/\partial z\}_1^k$ in
$\mathcal O(\Omega)$.
Thus, if $V$ is piecewise harmonic then
$\partial V/\partial z$
is piecewise holomorphic.
A  major result in this paper is the following converse: if $V$  is a 
subharmonic function
such that $\partial V/\partial z$
is piecewise holomorphic then $V$ is locally
piecewise harmonic.
More precisely, we prove:

\medskip

\noindent
{\bf Theorem 1.}
\emph{Let
$\partial V/\partial z $ be piecewise holomorphic 
with respect to a $k$-tuple $\{g_\nu\}_1^k$ in some open set
$\Omega$.
For each simply connected open subset
$U$ of $\Omega$ one can  choose  a $k$-tuple of harmonic functions
$H_1,\ldots,H_k$
such that $\partial H_i/\partial z=g_i$ in $U$, $1\le i\le k$. Moreover,
if $U_0$ is a relatively compact subset of $U$
there exists a finite number of constants 
$c_\nu=c_\nu(U_0)$, $1\leq \nu\leq m$, 
such that
the restriction $V|_{U_0}$ is piecewise harmonic with respect
to a subfamily of the  $m\cdot k$ many  harmonic functions
$\{ H_i+c_\nu\}$.}

\medskip

\noindent 
The proof of Theorem 1 requires several steps.
It is
based upon the results of \S 1 and \S 2 and will be completed 
only at the end of \S 3.
Let us point out that the
difficulty in proving Theorem 1 stems from the fact that no 
special assumption is
imposed on the open set $\Omega\setminus\text{supp}(\Delta(V))$, i.e., 
for a general null set
$K$ of $\Omega$ there may {\em \`a priori} exist a relatively compact
subset $U$ of $\Omega$ such that
the number of connected components of
$\Omega\setminus K$
which intersect $U$ is infinite.
The main burden in the proof of Theorem 1 is then to show
that this cannot occur when $V\in\text{SH}_0(\Omega)$ and
$K=\text{supp}(\Delta(V))$.

\medskip

\noindent
Our final topic is about algebraic functions. In \S 4 we make use of the previously developed material to 
prove a result about non-negative Riesz measures supported by compact null
sets in ${\bf{C}}$ whose Cauchy transforms satisfy an algebraic equation.
More precisely, let
$\mu$ be such a measure, denote by $K$ the support of $\mu$ and set
\[
\hat\mu(z)=\int\!\!\int_K\frac{d\mu(\zeta)}{z-\zeta}.
\]
We say that the Cauchy transform
$\hat\mu$ satisfies an algebraic equation
if there exist some $k\geq 1$ and polynomials
$p_0(z),\ldots,p_k(z)$ such that
\[
p_k(z)\cdot \hat\mu^k(z)+\ldots+p_1(z)\cdot \hat\mu(z)+p_0(z)=0, 
\quad z\in{\bf{C}}\setminus K.\tag{*}
\]
Note that {\em \`a priori} $K$ is just 
a null set and in general one can hardly say more than that. However, 
assuming (*) we can substantially improve this and get the following 
decription of $K$: 
\medskip

\noindent
{\bf Theorem 2.}
\emph{If (*) holds then the support of $\mu$ is a real analytic set of 
dimension at most one.}
\medskip

\noindent
Finally, in \S 5 we discuss some further directions, open problems and 
conjectures inspired by the topics treated in this paper.


\bigskip

\bigskip


\centerline{\bf\large 1. A Key Lemma}
\bigskip

\noindent
 Let $\Omega$ be an open  and connected set in
$\mathbf{C}$ and $V\in\text{SH}_0(\Omega)$. 
Set $K=\text{supp}(\Delta(V))$ and decompose $\Omega\setminus K=\bigcup_{\alpha\in C}\omega_\alpha$ into open connected components. Suppose that
$V=0$ in an open subset $U=\bigcup_{\alpha\in A}\omega_\alpha$ of $\Omega\setminus K$
and furthermore that
\[
\mathfrak{Re}(\partial V/\partial z)<0,\quad z
\in W:=\Omega\setminus(K\cup U)=\bigcup_ { \beta\notin A} \omega_\beta
\,.
\]
Observe that $\partial V/\partial z$ is a holomorphic function in each component $\omega_\alpha$.

\medskip
\noindent
{\bf 1.1.~Lemma.} {\em Let $z_0\in\omega_\alpha\subset U$ and assume that 
$\ell=\{z_0+s\,\colon\,0\leq s\leq s_0\}$ is a line segment contained in 
$\Omega$.
If $0<\delta<\text{{\em dist}}(\ell,\partial\Omega)$
and the open disk $D_\delta(z_0)$ of radius $\delta$ centered at
$z_0$ is contained in $\omega_\alpha$, then 
$$\{z\,\colon\, \text{{\em dist}}(z,\ell)<\delta\}=\bigcup_{0\leq s\leq s_0}D_\delta(z_0+s)\subset \omega_\alpha.$$}
\medskip

\noindent
{\bf Remark.} The subsequent proof uses methods similar to those of 
\cite[Lemma 2]{Be-Ru}, in particular the
idea to use the $\Psi$-function below. However, the new (and general) 
situation in Lemma 1.1  
is that no finiteness condition is imposed
on the range of $\partial V/\partial z$.  
\medskip

\noindent
\emph{Proof of Lemma 1.1.} By the choice of $\delta$, the set
$\{z\,\colon\, \text{dist}(z,\ell)<\delta\}$ is a 
relatively compact subset of $\Omega$. Let $\epsilon>0$ and define the holomorphic function
\[
\Psi(z)=\text{Log}(-\epsilon+\partial V/\partial z)\,,\quad\, 
z\in\Omega
\setminus K,
\]
\medskip
where the single-valued branch of the complex Log-function
is chosen so that
\[
\pi/2<\mathfrak{Im}\,\Psi<3\pi/2.
\]
This is clearly possible, since by assumption $-\epsilon+\partial V/\partial z\leq -\epsilon$ in $\Omega
\setminus K$. Furthermore, since $\partial V/\partial z$ is locally integrable, 
$\Psi\in L^1_{\text{loc}}(\Omega)$. 
\medskip

\noindent
Consider a non-negative cut-off  function $\rho$ supported by the unit disk with integral $1$. Let $\delta>0$, define $\rho_\delta(z)=\delta^{-2}\rho(z)$, and set 
$$
\Psi_\delta:=\text{Log}(-\epsilon+\rho_\delta *\partial V/\partial z).
$$
Taking a derivative, we get

\begin{equation}\label{1}
\begin{split}
&\partial\Psi_\delta/\partial\bar z=
\frac{1}{4}\cdot\frac{{\rho_\delta*\Delta V }}{-\epsilon+\rho_\delta *\partial V/\partial z}\implies\\
&\mathfrak{Re}(\partial\Psi_\delta/\partial\bar z)=
\frac{(-\epsilon+\rho_\delta *\mathfrak{Re}(\partial V/\partial z))\cdot \rho_\delta*\Delta V}{4
|\epsilon-\rho_\delta *\partial V/\partial z|^2}.
\end{split}
\end{equation}

\noindent
Since $\Delta V$ is a non-negative Riesz measure and  $\mathfrak{Re}(\partial V/\partial z)$ is a non-positive function, we deduce from (1) that $\mathfrak{Re}(\partial\Psi_\delta/\partial\bar z)$ is a non-positive function. Passing to the limit as $\delta\to 0$ we conclude that the distribution derivative 
$\mathfrak{Re}(\partial \Psi/\partial \bar z)$ is a non-positive Riesz measure.
Next, we can write
\begin{equation*}
\Psi(z)=\sigma(z)+i\tau(z),
\quad\, \pi/2<\tau(z)<3\pi/2,
\end{equation*}
where $\sigma(z)=\text{Log}|\epsilon-\partial V/\partial z|$ is the real part of $\Psi(z)$.

\medskip

\noindent
Let us now choose a non-negative
test function $\phi$ with compact support in
the disk $|z|\leq \delta$ such that $\phi(z)>0$ if $|z|<\delta$
and $\iint \phi(z)dxdy=1$.
By  the definition
of $\partial/ \partial\bar z$ the inequality $\mathfrak{Re}(\partial \Psi/\partial \bar z)\leq 0$ gives that
\begin{equation}\label{2}
\partial_x(\phi*\sigma)\leq\partial_y(\phi*\tau).
\end{equation}

\noindent

\noindent 

\noindent
Since $\pi/2\leq \tau\leq 3\pi/2$, the absolute value of the right-hand side 
is majorised by $M=\frac{3\pi}{2}\cdot||\partial_y(\phi)||_1$, where
$||\partial_y(\phi)||_1$ denotes the $L^1$-norm.
Next,  consider the function $s\mapsto \phi*\sigma(z_0+s)$,
where $0\leq s\leq s_0$. Applying (2) and setting $z_1=z_0+s_0$
 we obtain
\begin{equation}\label{4}
\frac{d}{ds}(\phi*\sigma(z_0+s))\leq M\implies
\phi*\sigma(z_1)\leq\phi*\sigma(z_0)+M\cdot s_0.
\end{equation}
Since $K=\text{supp}(\Delta(V))$ is a null set we can identify $\sigma$ with the
following $L^1_{\text{loc}}$-function 
\begin{equation}\label{3}
\sigma(z)=
\text{Log}|\epsilon|\cdot \chi_U+
\text{Log}|\epsilon-\partial V/\partial z|\cdot\chi_W,\quad
W=\Omega\setminus (K\cup U)\,.
\end{equation}

\noindent 
Set
$f_\epsilon
=\text{Log}\bigl|\epsilon-\partial V/\partial z\bigr|\cdot\chi_W
$.
From now on $\epsilon<1$ so that $\text{Log}|\epsilon | <0$. Since the support of $\phi$ is small enough (i.e., less than the distance $\delta $ from $z_0$ to the boundary) $\phi*\sigma(z_0)=\text{Log}|\epsilon |(\phi*\chi_U)(z_0)=\text{Log}|\epsilon |.$
Inserting in (4) the expression $f_\epsilon$, inequality (3) gives
\begin{equation}\label{5}
1\leq\phi*\chi_U(z_1)+
\frac{1}{\text{Log}|\frac{1}{\epsilon}|}\cdot[
-\phi*f_\epsilon(z_1)+M\cdot s_0]\,.
\end{equation}

\noindent
At this stage we perform a limit as $\epsilon\to 0$. 
For this note first that 
the function  $-\mathfrak{Re}(\partial V/\partial z)\cdot\chi_W$
belongs to  $L^1_{\text{loc}}$ and is $>0$ in $W$.
Moreover, the disk $D_\delta(z_1)$ is relatively compact in
$\Omega$.  
Elementary measure theory shows that for any 
$h\in L^1_\text{loc}(\Omega)$ such that $\mathfrak{Re}(h)\geq 0$ in $W$ and $\{\mathfrak{Re}(h)=0\}\cap W$
is a null set one has
\begin{equation}\label{6}
\lim_{\epsilon\to 0}\,
\frac{1}{\text{Log}\!\left|\frac{1}{\epsilon}\right|}\cdot 
\iint_{D_\delta(z_1)\cap W}\, |\text{Log}(|\epsilon+h|)|dxdy=0\,.
\end{equation}

\noindent
Apply this with $h=-\partial V/\partial z$. 
Since the
test function $\phi$ has support in $|z|\leq\delta$, we have
the inequality 
\begin{equation}\label{7}
|\phi*f_\epsilon(z_1)|\leq ||\phi||_\infty\cdot
\iint_{D_\delta(z_1)}\, |f_\epsilon(z)|dxdy.
\end{equation}

\noindent 
By (6) the quotient of this by ${\text{Log}\!\left|\frac{1}{\epsilon}\right|}$ tends to zero as
$\epsilon\to 0$.
So after a passage to the limit as $\epsilon\to 0$, it follows from (5) and (7) 
that
\begin{equation}\label{8}
1\leq \phi*\chi_U(z_1)\,.
\end{equation}

\noindent 
Finally, since $\phi(z)>0$ when $|z|<\delta$, inequality (8)
implies
that $D_\delta(z_1)\setminus U$ is a null set.
Hence the restriction of the
subharmonic function $V$  to this open disk is almost everywhere zero. 
Since subharmonic functions appear  as a subspace of $L^1_{\text{loc}}$-
functions
we conclude that $D_\delta(z_1)\subset U$. This completes the proof of 
Lemma 1.1.\hfill $\Box$

\medskip

\noindent
Lemma 1.1 suggests defining the following notion:
\medskip

\noindent
{\bf Definition.}
For every $z\in\Omega$ set
\[
\rho^*(z)=\max\{a\in (0,\infty)\,\colon\, z+t\in\Omega\text{ for all real } 
0<t<a\}.
\]
If 
$U$ is an open subset of $\Omega$ we define the 
\emph{forward star domain} of $U$ by 
\[ 
\mathfrak{s}^\uparrow(U)=\bigl\{ z\in\Omega\,\colon\,\exists\,\zeta\in U\,\,
\text{such that}\,\,z=\zeta+t\,\,\text{for some}\,\, 0\leq t<\rho(\zeta)\}.
\]
\medskip

\noindent
A more concise formulation of Lemma 1.1 is then as follows:

\medskip

\noindent
{\bf 
1.2.~Theorem.} \emph{Let $V\in\text{SH}_0(\Omega)$, $K=\text{supp}(\Delta(V))$
and assume that $\Omega\setminus K$ is the disjoint union $U\cup W$ of two open sets such that
$\mathfrak{Re}(\partial V/\partial z)<0$ in
$W$ and $V=0$ in $U$.
Then $U=\mathfrak{s}^\uparrow(U)$.}

\medskip

\noindent Notice that Theorem 1.2 applies to an arbitrary subharmonic function in $SH_0(\Omega)$, not necessarily piecewise harmonic. It will be crucial for our study of the piecewise
holomorphic case in \S 3 as well as for our next result that we proceed to 
describe.

\medskip

\noindent
{\bf 1.3.~Local subharmonic configurations.}~Let $V\in\text{SH}_0(D)$ and assume that 
$\partial V/\partial z$ is piecewise holomorphic with respect to some $k$-tuple
$g_1,\ldots,g_k$ in
$\mathcal O(D)$, where $D$ is an open disk centered
at the origin. With
$K=\text{supp}(\Delta(V))$ we further define the open subset
$U_\nu$ of $D\setminus K$
as the union of those connected components of
$D\setminus K$ where
$\partial V/\partial z=g_\nu$.
We assume that the origin belongs to the closure of
every $U_\nu$.
In the simply connected
disc $D$ we choose the unique $k$-tuple of harmonic functions
$H_1,\ldots,H_k$ satisfying 
\[
\partial H_\nu/\partial z=g_\nu\text{ and } H_\nu(0)=0,\quad
1\leq\nu\leq k.
\]
\noindent
Next, consider the $k$-tuple
$(g_1(0),\ldots,g_k(0))$ and the convex set $P$ generated by these complex 
numbers. Assume that
$g_k(0)$ is an extreme point of $P$.
This gives some
$\theta_*$ such that
\[
\mathfrak{Re}\!\left(e^{i\theta_*}\cdot g_\nu(0)\right)<
\mathfrak{Re}\!\left(e^{i\theta_*}\cdot g_k(0)\right),\quad
1\leq \nu\leq k-1.
\]
After a rotation if necessary we may further assume that
$\theta_*=0$ and thus (by continuity) there exists 
$\delta>0$ such that
\[
\mathfrak{Re}\!\left(e^{i\theta}\cdot g_\nu(0)\right)<
\mathfrak{Re}\!\left(e^{i\theta}\cdot g_k(0)\right),\quad 
1\leq \nu\leq k-1,\,\,-\delta<\theta<\delta.
\]
We can apply Theorem 1.2  to the subharmonic function
$e^{i\theta}\cdot(V-H_k)$ for $-\delta<\theta<\delta$, and setting $U=\{V=H_k\}$ we conclude:
\medskip

\noindent
{\bf 1.4 Proposition.}
\emph{If $\square=\{(x,y)\in\mathbf{R}^2\,\colon\, -a<x,y<a\}$ and $a>0$ is 
sufficiently small then the domain 
$\square\cap U$ is connected
and given by}
\[
\square\cap U=\{(x,y)\in\square\,\colon\, x>\rho(y)\},
\]
\emph{where $\rho(0)=0$ and $\rho$ is a 
 Lipschitz continuous function of norm
$\leq \frac{\text{cos}\,\delta}{{\text{sin}\,\delta}}$.}

\medskip

\noindent
A similar conclusion holds for
other indices as well. Indeed, if
$g_\nu(0)$
is an extreme point of $P$ for every $1\leq\nu\leq k$ then
we obtain open connected sets
$U_1,\ldots,U_k$
as above after suitable rotations. This leads to the following result.

\medskip

\noindent
{\bf 1.5.~Theorem.}
\emph{Let $V\in\text{SH}_0(D)$ and assume that 
$\partial V/\partial z$ is piecewise holomorphic with respect to some $k$-tuple
$g_1,\ldots,g_k$ in
$\mathcal O(D)$, where $D$ is an open disk centered
at the origin. Assume further that each $g_i(0)$ is an extreme point of the convex hull $P$ of $(g_1(0),\ldots,g_k(0))$. Then there exists $c\in \bf R$ such that  in a neighborhood of the origin one has 
$V=\text{{\em max}}(H_1,\ldots, H_k)+c$.}
\medskip

\noindent \emph{Proof.}~The hypothesis implies that for each given 
$1\leq m\leq k$
there exists some
$\theta$ such that
$\mathfrak{Re}(e^{i\theta}g_\nu)<\ldots<
\mathfrak{Re}(e^{i\theta}g_m)$, 
$\nu\neq m$. 
Theorem 1.2 applies after a rotation. It follows that $U_m\cap D(\delta)$
is  connected for a sufficiently small $\delta$. Since this holds for 
every $m$  it
follows that $V$ is piecewise harmonic with respect to the $k$-tuple
$H_1,\ldots,H_k$
in  $D(\delta)$.
There remains to see that $V$ is the maximum function. For this we may consider
without loss of generality the index $m=1$. After a rotation we find that there exists a function
$\rho(y)$ such that
\[ U_1=\{ (x,y)\in D(\delta)\,\colon\,x>\rho(y)\}\,\text{ and }\, 
\partial_x H_\nu <\partial_x H_1,\,\nu\geq 2\,.
\]
\noindent
We have  to show that
$H_1(x,y)<V(x,y)$ when $x<\rho(y)$.
To do this we fix $y_0$ and consider the function
$
x\mapsto V(x,y_0)$.
When $x<\rho(y_0)$ the partial derivative
$\partial_x(V)$ is equal to $\partial_x(H_\nu)$ for some
$\nu\geq 2$
on intervals
outside some finite set
where $V$ may shift from one $H$-function to another when 
a level curve
$\{H_i=H_\ell\}$ intersects the line $y=y_0$.
By the strict inequalities above 
$x\mapsto V(x,y_0)-H_1(x,y_0)$
is strictly decreasing and since it is zero when
$x=\rho(y_0)$ Theorem 1.5 follows.\hfill $\Box$

\medskip

\noindent
{\bf 1.6.~A relaxed assumption.}
Let us drop the hypothesis
that the origin belongs to
$\bar U_\nu$ for every $\nu$
and suppose instead that
there is some $1\leq \ell\leq k-1$ such that
the extreme points of $P$ are
$g_i(0)$, $1\leq i\leq\ell$.
Without loss of generality we may assume that
the origin belongs to $\bar U_i$, 
$1\leq i\leq\ell$, and that the vertices of $P$ are 
labelled consecutively $g_1(0),\ldots,g_\ell(0)$ in say 
counter-clockwise order.
The example in
\S 2.11 below shows that in this case we cannot conclude that
$V$ is given by
the maximum of $H_1,\ldots,H_\ell$ up to a constant.
However, the following extension of Theorem 1.5 holds:
\medskip

\noindent
{\bf 1.7 Theorem.}
\emph{Suppose as above that
$\{g_i(0)\}_1^\ell$
are the extreme points of $P$ and that for $i\in \{\ell+1,\ldots,k\}$ one has
$$
g_i(0)\notin \bigcup_{j=1}^{\ell}
\big\{(1-\alpha)g_{[j]}(0)+\alpha g_{[j+1]}(0):
0\le \alpha\le 1\big\},
$$
where $[j]=j$ for $1\le j\le \ell$ and $[\ell+1]=1$. Then in a sufficiently 
small neighborhood of the origin one has 
$V=\text{{\em max}}(H_1,\ldots,H_\ell)$ 
up to a constant.}


\bigskip

\bigskip

\centerline{\bf\large 2.~Subharmonic Configurations:} 
\centerline{\bf\large The General Piecewise Harmonic Case}
\bigskip

\noindent
We begin with some preliminary observations which follow
from the maximum principle for  subharmonic functions and Stokes' Theorem.
We then study harmonic level sets and give a local description of  
arbitrary subharmonic configurations.
\medskip

\noindent
Let $H_1,\ldots,H_k$ be harmonic functions in $\Omega$ and 
$V\in\text{SH}_0(\Omega)$
be piecewise harmonic function with respect 
to this $k$-tuple.
In $\Omega$ we get the real analytic set
\[ 
\Gamma=\bigcup_{i\neq j}\,\{H_i=H_j\}.
\]
\noindent
Let $\{U_\alpha\}$ be the connected components of
$\Omega\setminus\Gamma$. Then we have:
\medskip

\noindent 
{\bf 2.1.~Lemma.}
\emph{For each $\alpha$ there exists $1\leq i(\alpha)\leq k$ such that
$V=H_{i(\alpha)}$ in $U_\alpha$.}
\medskip

\noindent 
\emph{Proof.}~Given $U_\alpha$ there is some permutation of the 
indices such that
\[ 
H_{j(1)}<\ldots <H_{j(k)}.
\]
Set $K=\text{supp}(\Delta(V))$. For 
each $1\leq i\leq k$
we define 
\[
U_\alpha(i)=\{ z\in U_\alpha\setminus K\,\colon\, 
V=H_i\text{ in some neighborhood of}\,\, z\}. 
\]
By assumption one has 
$\bigcup_i\, U_\alpha(i)=U_\alpha\setminus K$.
Let
$m$ be the largest integer such that $U_\alpha(j(m))$
is non-empty. Then we have:
\medskip

\noindent
\emph{Sublemma.}~\emph{The set $U_\alpha(j(m))$ is dense in
$U_\alpha$}.
\medskip

\noindent
\emph{Proof.}~Assume the contrary and set
$U_*=U_\alpha\setminus \bar U_\alpha(j(m))$.
Since $U_\alpha$ is connected we cannot have
$U_\alpha(j(m))\cup U_*=U_\alpha$ and hence there exists a point
\[ 
p_*\in\partial (U_\alpha(j(m)))\cap\, U_\alpha\,.
\]
Consider a point 
$p\in U_\alpha(j(m))\cap D$ very close to $p_*$.
Let $D$ be a disc centered at $p$ of some radius $r$
whose closure stays in
$U_\alpha$. With  $p$ sufficiently  close to $p_*$
the set $D\cap U_*$ is non-empty.
The mean value inequality for subharmonic functions gives
\[
H_{j(m)}(p)= V(p)\leq\frac{1}{\pi r^2}\iint_D\,
V(x,y)\cdot dxdy\,.\tag{9}
\]
Set $H_*(z)=\max(H_{j(1)},\ldots,H_{j(m)-1})$.
Since $D\cap U_*$ is non-empty and 
$K$ is a null set we have
$
V(z)\leq H_*(z)
$ almost everywhere in
$D\cap U_*$.
But then (9) cannot hold since we have the strict inequality
$H_*<H_{j(m)}$.\hfill $\Box$

\medskip

\noindent
\emph{Proof of Lemma 2.1, continued.}~By the 
Sublemma $U_\alpha(j(m))$ is dense in $U_\alpha$
and since all the sets $U_\alpha(i)$ are open
we have
\[ 
U_\alpha(j(m))= U_\alpha\setminus K\,.
\]
This means that the
$L^1_{\text{loc}}$-function $V$ equals $H_{j(m)}$
in the whole set $U_\alpha$ and then Lemma~2.1 follows with
$i(\alpha)=j(m)$.\hfill $\Box$

\medskip

\noindent
{\bf Remark.}~Note that Lemma 2.1 gives the inclusion
\[
\text{supp}(\Delta(V))\subset \Gamma.
\]
\noindent
Another way of proving Lemma 2.1 is by means of 
Grishin's Lemma \cite{Grishin}, see also \cite{Fuglede}.
In fact, using \cite{Grishin} one can 
show that if $V\in \text{SH}_0(\Omega)$ is 
piecewise harmonic then $\text{supp}(\Delta(V))$ is a null set, so the latter 
property need not be assumed already from the start (which we did for 
the reader's convenience).
\medskip

\noindent
{\bf 2.2.~A description of
$\Delta(V)$.}~Consider some pair $(U_\alpha,U_\beta)$
with $i(\alpha)\neq i(\beta)$
and such that $\partial U_\alpha\cap\,\partial U_\beta\neq \emptyset$.
As explained in \S 2.5 below,
the portion of this common boundary set that avoids the closed union
of the remaining $U$-sets is a smooth real analytic curve
$\gamma$ possibly up to a discrete set.
Let $ds_\gamma$ be  arc-length measure on
$\gamma$ and suppose 
$H_{i(\alpha)}>H_{i(\beta)}$ holds in $U_\alpha$
while $H_{i(\alpha)}<H_{i(\beta)}$ in $U_\beta$.
Along $\gamma$ we choose the normal $\mathfrak{n}_\gamma$
directed into $U_\alpha$.
Hence  the normal derivatives satisfy
\[
\partial_{\mathfrak{n}_\gamma}H_{i(\alpha)}>0\,\,\text{and}\,\, 
\partial_{\mathfrak{n}_\gamma}H_{i(\beta)}<0
\]
outside the discrete
set of possible singularities for the
level curve $\{H_{i(\alpha)}=H_{i(\beta)}\}$.
With these notations Stokes' Theorem gives:
\medskip

\noindent {\bf 2.3.~Proposition.}~\emph{One has}
$
\Delta(V)|_\gamma=
\bigl [\partial_{\mathfrak{n}_\gamma}H_{i(\alpha)}-
\partial_{\mathfrak{n}_\gamma}H_{i(\beta)}\bigr]\cdot ds_\gamma$.
\medskip

\noindent {\bf 2.4.~Remark.}~Let $G,H$ be a pair of harmonic functions
defined in some domain $\Omega$, set 
$\Gamma=\{G=H\}$ and let $p\in\Gamma$ be a regular point,  
i.e., $\nabla(G)(p)-\nabla H(p)\neq 0$.
Consider a small disk $D$ centered at $p$ and the two domains
\[
U_+=\{G>H\}\text{ and } U_-=\{G<H\}.
\]
Then $V=\text{max}(G,H)$ is subharmonic while
the opposed function $\text{min}(G,H)$ fails to be
subharmonic. The lesson of this observation is that
when the pair $G,H$ appears in a  configuration
of a subharmonic function $V$  their normal  derivatives satisfy 
\[
\partial_{\mathfrak{n}}G\geq
\partial_{\mathfrak{n}}H, 
\] 
where $\mathfrak{n}$ is the normal to $\Gamma$ \emph{directed into}
$U_+$.
This simple -- but essential -- observation
will be frequently used later on.

\medskip

\noindent {\bf 2.5.~Harmonic level sets.}~Let 
$H(x,y)$ be a harmonic function defined in some open disk $D$ centered at the 
origin in $\mathbf{C}$ and $z=x+iy$ be the complex variable.
Now $H=\mathfrak{Re}(g)$ for some $g\in\mathcal O(D)$.
If $g$ vanishes of some order $m\geq 1$ at $z=0$ there exists a conformal map 
$\rho(\zeta)$ from a disk in the complex $\zeta$-plane such that
$g\circ\rho(\zeta)=\zeta^m$.
The zero set of $\mathfrak{Re}(\zeta^m)$ is the union of 
lines $\text{arg}(\zeta)=\frac{\pi}{2}+\nu\mathfrak{\pi}{m}$, 
$0\leq\nu\leq m-1$. Passing to the $z$-disk
and shrinking $D$ if necessary 
we get that 
$\{H=0\}$ is the union of $m$ smooth real analytic curves 
$\gamma_1,\ldots,\gamma_m$ and  $D\setminus \{H=0\}$ consists of $2m$ 
pairwise disjoint
open sets $U_1,\ldots\,U_{2m}$, each $U_\nu$ being 
bordered by a pair of $\gamma$-curves intersecting at the origin where 
the angle between their tangential vectors is $\frac{\pi}{m}$.
Thus, every $U_\nu$ is  a simply connected real analytic sector.

\medskip

\noindent
Let us now consider a finite family of (distinct) harmonic 
functions $H_1,\ldots,H_k$ in $D$ satisfying $H_\nu(0)=0$ for all $\nu$.
Set
\begin{equation*}
\Gamma=\bigcup_{i\neq \nu}\,\{H_i-H_\nu=0\}.
\end{equation*}

\noindent
Applying the previous observation to all pairs $(H_i,H_\nu)$ it
follows that $\Gamma$ is a finite union of smooth 
real analytic curves $\gamma_1,\ldots,\gamma_M$ such that 
they all pass through the origin and are
pairwise disjoint in the punctured disk 
$$
\dot{D}=\setminus\{(0,0)\}.
$$
Of course, in general one must shrink $D$ 
to achieve this. Thus, provided that $D$ is sufficiently small, 
$D\setminus\Gamma$ is a union of pairwise disjoint real analytic sectors, each
of which is
bordered by two "half-curves" coming from the above family of $\gamma$-curves. 
\medskip

\noindent
Notice that no  special assumptions are imposed on the gradient vectors
of the $H$-functions at the origin. For example,
they may all  be zero. It
may therefore occur that some of the real analytic sectors $\Omega$
are bordered by a pair of $\gamma$-curves which 
\emph{do not} intersect transversally at the origin. Up to a conformal map a 
typical topological picture is that a real analytic sector is given by
\begin{equation*}
\Omega=\{(x,y)\,\colon\,0<x<\delta,\,0<y<\rho(x)\},
\end{equation*}
where $\rho(x)$ is a positive real analytic function on $(0,\delta)$
and there exists a holomorphic function $g$ in $D$ such that
$\mathfrak{Re}(g(x,\rho(x)))=0$. 
\medskip

\noindent
{\bf 2.6.~Local subharmonic configurations.}~Given an open disk $D$ centered 
at the origin and a $k$-tuple of harmonic functions 
$H_1,\ldots,H_k$ as above, we consider some $V\in\text{SH}_0(D)$ which is 
piecewise harmonic with respect to this  $k$-tuple.
Lemma 2.1 implies that $\text{supp}(\Delta(V))$
is contained in the set $\Gamma$ defined at the beginning of this section.
Hence, if $\omega_1,\ldots,\omega_N$ are the real analytic sectors whose 
union is $D\setminus\Gamma$, we find for each $\omega_\nu$ some 
$1\leq j(\nu)\leq k$ such that $V=H_{j(\nu)}$ in $\omega_\nu$.
\medskip

\noindent
Next we describe the positive measure $\Delta(V)$.
Outside the origin it is supported by (a subset of) $\Gamma$ and 
Proposition 2.3 shows
that if one has two adjacent $\omega$-sectors, say $\omega_1,\omega_2$ with 
$j(1)\neq j(2)$, then the portion of $\Delta(V)$ supported by the 
real 
analytic curve
$\gamma=\partial \omega_1\cup\partial\omega_2$ is 
the positive measure 
\[
[\partial_{\mathfrak{n}_\gamma}H_{j(1)}-\partial_{\mathfrak{n}_\gamma}
H_{j(2)}]\cdot ds_\gamma,
\]
where $ds_\gamma$ 
is arc-length measure and $\mathfrak{n}_\gamma$ is the normal to $\gamma$  
directed into $\omega_1$ when $H_{j(1)}>H_{j(2)}$ holds in $\omega_1$
while
$\mathfrak{n}_\gamma$ changes sign and is directed
into $\omega_2$ if it happens that
$H_{j(2)}>H_{j(1)}$ holds in $\omega_1$, see Remark 2.4.
There remains to show that $\Delta(V)$ cannot contain a point mass at the 
origin. 
For this, we construct
the logarithmic potential $W$ of $\mu=
\Delta(V)|_{\dot{D}}$.  Note that since $ \Delta(V)|_{\dot{D}}$  is a locally real-analytic density on real-analytic curves $W$ is a continuous and bounded
subharmonic function and $V-W$ is harmonic outside the origin. So if 
$\Delta(V)$ has a point mass at the origin, there exists a constant $a>0$ 
such that 
$V=a\text{Log}(|z|)+W+G$, where $G$ is harmonic in $D$.
This is impossible since $V$ is a bounded function in the punctured 
open disk $\dot{D}$. 
We conclude that $V$ can be taken as a continuous function, i.e., 
we have proved:

\medskip

\noindent {\bf 2.7.~Theorem.}~\emph{Every 
piecewise harmonic subharmonic function is continuous}.
\medskip
\newline
\noindent
Let us summarize our results so far. Given (distinct) 
harmonic functions 
$\{H_i\}_{1}^k$ and a subharmonic function $V$ which is piecewise
harmonic with respect to this family, the following holds
if the disk $D$ (centered at the origin) is sufficiently small:
\medskip

\noindent {\bf 2.8.~Theorem.}~\emph{There exists a finite family of
disjoint
real analytic sectors, say $\omega_1,\ldots,\omega_m$ such that for 
each $1\leq i\leq m$ one has}
\[
V|_{\omega_i}=H_{j(i)},\quad 1\leq j(i)\leq k\,.
\]
\emph{Moreover, when $1\leq i\leq m-1$ a half-arc  $\gamma_i$
from the 
level set $\{H_{j(i+1)}=H_{j(i)}\}$
borders $\bar\omega_{i+1}\cap\bar\omega_i$ outside the origin 
and here one has the strict inequality}
\[\partial_{\mathbf{n}_i} H_{j(i+1)}>
\partial_{\mathbf{n}_i} H_{j(i)},
\]
\emph{where $\mathbf{n}_i$ is the normal to $\gamma_i$ directed into 
$\omega_{i+1}$.
When $i=m$ one returns from $\omega_m$ to $\omega_1$ and here one has 
$\partial_{\mathbf{n}_m}  H_{j(1)}>
\partial_{\mathbf{n}_m} H_{j(m)}$, where  $\mathbf{n}_m$ is the 
normal to a half-arc
$\gamma_m$  of the level set $\{H_{j(m)}=H_{j(1)}\}$
which is  directed
into $\omega_1$. Finally, the measure $\Delta(V)$ is given by}
\[\Delta(V)=\sum_{i=1}^{m-1}\, [\partial_{\mathbf{n}_i} H_{j(i+1)}-
\partial_{\mathbf{n}_i} H_{j(i)}]\cdot d_{\gamma_i}s+
 [\partial_{\mathbf{n}_m} H_{j(1)}-
\partial_{\mathbf{n}_m} H_{j(m)}]\cdot d_{\gamma_1}s\,.
\]

\noindent
{\bf 2.9.~A non-transversal  case.}~Let 
$H_1,\ldots,H_k$ be harmonic in an open 
disk $D$ centered at the origin.
Assume that they are all zero at the origin and
their gradients there satisfy
\[ 
\nabla(H_\nu)=(0,b_\nu),\quad b_1<\ldots<b_k.
\]
In this case the above results give 
a transparent description of all subharmonic configurations
with respect to this $k$-tuple in a small neighborhood of the origin.
Indeed, let $V$ be such a subharmonic
configuration.
Assume that the closure of the two sets
$U_1=\{V=H_1\}$ and $U_k=\{V=H_k\}$
both contain the origin.
Theorem~2.8 implies that $U_1$ contains
a sector of the form 
\[
\Omega_+=\{(x,y)\colon\,\, 0<x<\delta,\,|y|<ax\}
\] 
for some appropriate $a,\delta>0$. Similarly, 
$U_k$ contains a sector $\Omega_-$ where $x<0$ and
$y<a|x|$.
In the upper semi-disk $D_+$ where $y>0$
we have smooth half-arcs
\[
\gamma^+_{i\nu}=\{(x,y)\colon\,\, H_i(x,y)=H_\nu(x,y),\, y>0\}
\]
and similar half-arcs $\gamma_{i\nu}^-$ in the lower semi-disk $D_-$.
With these notations we have:

\medskip

\noindent {\bf 2.10.~Theorem.}~\emph{Let 
$V$ be a subharmonic configuration such that
$\bar U_1$ and $\bar U_k$ contain $(0,0)$. There 
exist  integers $m,n\geq 2$, 
a pair of sequences $1= j^+_1<\ldots <j^+_m=k$, 
$1=j_1^-<\ldots<j_n^-=k$, and some $\delta>0$ such that}
\[
V|_{D_+(\delta)}=\text{max}( H_{j_1^+},\ldots,H_{j_m^+})\text{ and }
V|_{D_-(\delta)}=\text{max}( H_{j_1^-},\ldots,H_{j_n^-}).
\]
\emph{Conversely, every such pair of $j$-sequences yields 
a subharmonic configuration.}

\medskip

\noindent
{\bf 2.11.~Example.}~Let $k=3$, $H_{1}(x,y)=0$, 
$H_{2}(x,y)= 4x+x^2-y^2$, $H_{3}(x,y)=-x$. There are 
three level curves through $(0,0)$ to functions of the form $H_i-H_j$ with 
$i\neq j$. These are depicted in Figure 1 below. 

\begin{figure}[!htb]
\begin{center}
\includegraphics[scale=0.75]{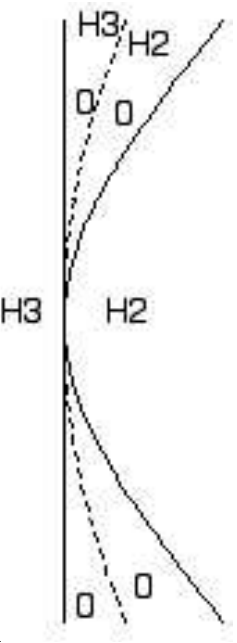}
\end{center}
\caption{Maximal and non-maximal subharmonic configurations.}
\end{figure}

\noindent
Here we get three different 
subharmonic configurations (when
the origin also belongs to the closure of $\{V=H_2\}$), one of these 
configurations being $\text{max}(H_1,H_2,H_3)$. 
The function in the figure closest to the origin in each sector is the 
restriction of $V$ to that sector.


\bigskip

\bigskip

\centerline{\bf\large
3.~Piecewise Holomorphic Functions}
\bigskip

\noindent
It suffices to prove Theorem 1 locally, i.e., we can restrict our 
attention to an open neighborhood of the origin where
$\partial V/\partial z$ is piecewise holomorphic with respect to some $k$-tuple 
$g_1,\ldots, g_k$. The neighborhood in question is chosen as an open
square
\[ 
\square=\{(x,y):\,\,-a<x,y<a\}.
\]
We shall first consider the case when
the following holds in
$\square$:
\[ 
\mathfrak{Re}(g_1)<\ldots<\mathfrak{Re}(g_k).\tag{**}
\]
\smallskip
Given a $k$-tuple of constants
$c_1,\ldots,c_k$
there exist 
harmonic functions $H_1,\ldots,H_k$ in $\square$ such that
$\partial H_\nu/\partial z=g_\nu$ and $H_\nu(0)=c_\nu$, $1\le \nu\le k$.

\medskip

\noindent
{\bf 3.1 Proposition.}~\emph{If (**) holds there exists
an open disk $D$ centered at $(0,0)$ 
such that $V|_D$ is piecewise harmonic with respect to
$H_1,\ldots, H_k$ up to additive constants.}

\medskip

\noindent
The proof requires several steps.
Set $K=\text{supp}(\Delta(V))$ and let $U_k$
be the subset of $\square\setminus K$ where
$\partial V/\partial z=g_k$.
Without loss of generality we may assume that all $g$-functions are active
in the sense that  
the sets  $\{\partial V/\partial z=g_\nu\}$ contain points arbitrarily close to
the origin for each $\nu$.
Theorem 1.2 applied to the subharmonic function $V-H_k$
gives $U_k=\mathfrak{s}^\uparrow(U_k)$, where $\mathfrak{s}^\uparrow(U_k)$ is 
the forward star domain of $U_k$ as 
defined in Section 1.
Since the origin by assumption   belongs to $\bar U_k$
this equality yields
\[ 
U_k=\{(x,y)\in\square:\,\,x>\rho^*(y)\},
\]
where $\rho^*(0)=0$. Moreover, as explained in Proposition 1.4,  $\rho^*$
is
Lipschitz continuous if we from the start shrink $\square$ a bit so that 
the inequality
\[
\mathfrak{Re}(e^{i\theta} g_\nu)<
\mathfrak{Re}(e^{i\theta} g_k),\quad -\theta_0<\theta<\theta_0,
\] 
holds in $\square$ for some $\theta_0>0$ and every $1\leq\nu\leq k-1$.
\medskip

\noindent
Since $U_k$ is connected $V-H_k$ is 
a constant function in $U_k$. We choose $c_k$
above so that $V=H_k$ holds in $U_k$.
Reversing signs in the Key Lemma 1.1 and considering the
subharmonic function $H_1-V$ it follows that if we set 
$U_1=\{\partial V/\partial z=g_1\}$ then
\[ 
U_1=
\{(x,y)\in\square:\,\,x<\rho_*(y)\},
\]
where $\rho_*$ is also Lipschitz continuous.
Since $U_1$ is connected $V-H_1$ is 
a constant function in $U_1$ and we choose $c_1$
such that $V=H_1$ holds in $U_1$.
\medskip

\noindent 
To complete the proof of Proposition
3.1 we  proceed by induction over $k$. Consider  first
the case
$k=2$. Then $\rho^*=\rho_*$
and by  Lipschitz continuity
the curve $x=\rho^*(y)$ is a null set. We
conclude that 
$V=\text{max}(H_1,H_2)$ in a small disk centered at the origin, as required.
\bigskip

\noindent {\bf The case $k\geq 3$.}~Since  
we only assume that $K=\text{supp}(\Delta(V))$ is a 
null set it is {\em \`a priori} not clear why the open subsets
$U_\nu$ of $\square\setminus K$
where $\partial V/\partial z=g_\nu$ have a finite number of
connected components when $2\leq\nu\leq k-1$.
To prove this
we shall consider the real analytic curve
$\Gamma=\{H_1=H_k\}$.
Since $\partial_x H_1<\partial_x H_k$ this curve is defined
by an equation of the form 
$x=\rho(y)$,
where $\rho$ is real analytic. Moreover, it is obvious that
\[ 
\rho_*(y)\leq\rho(y)\leq\rho^*(y)
\]
in some sufficiently small interval $-y_0<y<y_0$.
\medskip

\noindent
The $\Gamma$-curve is oriented by increasing $y$. 
The tangential derivatives $\partial_\Gamma H_\nu$
are real analytic functions on $\Gamma$ for each $\nu$.
We shall first consider these
tangential derivatives along the portion of $\Gamma$ where $y>0$.
Since the zero set of a real analytic function is discrete, it follows that
if $\square$ is if necessary decreased a bit (i.e., for $a$ small enough) 
then for every
$2\leq\nu\leq k-1$ the
function
\[ 
y\mapsto\partial_\Gamma H_k(\rho(y),y)-\partial_\Gamma H_\nu
(\rho(y),y)\,,\quad0<y<a,\tag{i}
\] 
is either  identically zero or else
strictly monotone, i.e., strictly increasing or decreasing.
Similarly, there exists a permutation of $\{2,\ldots,k-1\}$ such that
\[
\partial_\Gamma H_{j(2)}(\rho(y),y)\leq
\ldots\leq \partial_\Gamma H_{j(k-1)}(\rho(y),y)\,,\quad
0<y<a\,,\tag{ii}
\]
where $\leq$ means that we either have equality on
the whole portion of $\Gamma$ or a strict inequality.
\medskip

\noindent
{\bf 3.2.~The Non-Return Lemma.}~If there exists some $\delta>0$
such that $\rho_*(y)=\rho^*(y)$
holds for $0\leq y<\delta$ then  $V$ restricted to the
rectangle $\{-a<x<a,\, 0<y<\delta\}$ is equal to
$\text{max}(H_k,H_1)$ and we are done.
Next, we consider the situation when
no such $\delta$ 
exists.

\medskip

\noindent 
{\bf 3.3.~Lemma.}~\emph{Assume that $\rho^*-\rho_*$
is not identically zero on some interval $[0,\delta)$, that
strict inequalities hold in
(ii) and that the function in (i) is strictly monotone for some
$\delta>0$.
Then there exists $0<\delta_0<\delta$ such that}
$\rho_*(y)<\rho^*(y)$ for $0<y<\delta_0$.

\medskip

\noindent \emph{Proof.}~If 
the assertion is not true there exists a sequence of disjoint intervals
$J_\nu=(\alpha_*(\nu),\alpha^*(\nu))$
on the positive $y$-axis
which decrease to $y=0$ as $\nu\to\infty$ and such that at the end-points 
one has 
\[ 
\rho_*(\alpha_*(\nu))=
\rho^*(\alpha_*(\nu))\,\,\text{\and}\,\,
\rho_*(\alpha^*(\nu))=
\rho^*(\alpha^*(\nu))
\]
while $\rho_*(y)<\rho^*(y)$
holds inside every $J$-interval.
For any $\nu$ we consider the domain
\[ 
\Omega_\nu=\{(x,y)\,\colon\,  \rho_*(y)<x<\rho^*(y)\,\text{and}\,
\alpha_*(\nu)<y<\alpha^*(\nu)\}
\]
Inside each $\Omega_\nu$ we notice that $\partial V/\partial z$ is 
piecewise holomorphic with respect to
$g_2,\ldots,g_k$. By the induction assumption 
we may assume that $V$  is locally piecewise harmonic
in $\Omega_\nu$ with respect to $H_2,\ldots, H_{k-1}$ up to additive constants.
Therefore, when $\alpha_*(\nu)<y<\alpha^*(\nu)$
is kept fixed  the function
$
x\mapsto V(x,y)
$ 
is piecewise real analytic and  $\partial_x V$ is equal to
some $\partial_x H_\nu$ 
with $2\leq\nu\leq k-1$ 
outside a discrete set.
Since $\partial_x H_\nu<\partial_x H_1$ for each such $\nu$, we conclude that
the function
\[
x\mapsto H_1(x,y)-V(x,y)
\]
is strictly increasing. In the same way we find that
\[
x\mapsto H_k(x,y)-V(x,y)
\]
is strictly decreasing. Using this we conclude that the $\Gamma$-curve passes 
through $\Omega_\nu$, i.e., we must have
\[ 
\rho_*(y)<\rho(y)<\rho^*(y)\,\,\text{and}\,\, 
\alpha_*(\nu)<y<\alpha^*(\nu)\,.
\]
\medskip
\noindent
Here
$\rho^*=\rho=\rho_*$ at the end-points of $J_\nu$.
Set 
\[
p_\nu=(\rho(\alpha_*(\nu)),\alpha_*(\nu))\,\,\text{and}\,\, 
q_\nu=(\rho(\alpha^*(\nu)),\alpha^*(\nu))\,.
\]

\noindent
\emph{Sublemma.}~\emph{There cannot exist
a single index $2\leq j\leq k-1$ and a constant
$c_j$ such that
$V=H_j+c_j$ holds along $\Omega_\nu\cap\,\Gamma$.}
\medskip

\noindent
\emph{Proof}.~If this occurs we 
get a contradiction as follows.
First, by the induction over $k$ the restriction of
$V$ to $\Omega_\nu$ is piecewise harmonic
which yields a uniform bound for
$\partial_x(V)$.
Moreover, when
$\alpha_*(\nu)<y<\alpha^*(\nu)$
we encounter the point
$p=(\rho^*(y),y)$ where there exists an open neighborhood
such that
$\partial V/\partial z$ is piecewise holomorphic with respect to
$g_2,\ldots,g_k$. So by an induction over
$k$ it follows that
$V$ is piecewise harmonic in a neighborhood of
this point hence also continuous by Theorem 2.7.
Next, from the uniform bound of
$\partial_x V$ we have a constant $C$ which can be taken as the
maximum over sup-norms of $\{\partial_x H_\nu\}_2^{k-1}$
in a fixed neighborhood of the origin and
get
\[ 
\bigl|V(\rho(y),y)-V(\rho^*(y),y)\bigr|\leq C\cdot |\rho^*(y)-\rho(y)|
\quad\text{for all}\,\, \alpha_*(\nu)<y<\alpha^*(\nu)\,.
\]
Here $V(\rho^*(y),y)= H_k(\rho^*(y),y)$. So if
$V=H_j+c_j$ is valid on $\Omega_\nu\cap\Gamma$
for some constant $c_j$
we obtain
\[
\bigl|H_j(\rho(y),y)+c_j-H_k(\rho^*(y),y)\bigr|\leq C\cdot |\rho^*(y)-\rho(y)| 
\text{ whenever }\alpha_*(\nu)<y<\alpha^*(\nu).
\]
Passing to the limit as $y\to \alpha^*(\nu)$
or $y\to \alpha_*(\nu)$ we conclude that one has the two equalities:
\begin{equation*}
\begin{split}
&H_j(\rho(\alpha_*(\nu)), \alpha_*(\nu))+c_j=
H_k(\rho(\alpha_*(\nu)), \alpha_*(\nu)),\\
&H_j(\rho(\alpha^*(\nu)), \alpha^*(\nu))+c_j=
H_k(\rho(\alpha^*(\nu)), \alpha^*(\nu)).
\end{split}
\end{equation*}

\noindent
These two identities cannot hold if
$H_k-H_j$ is strictly monotone along
$\Gamma$. So there remains only the possibility that
$H_k-H_j$ is constant along $\Gamma$.
But this again gives a contradiction. For then we get
\[ 
V(z)=H_j(z)+c_j=H_k(z),\quad z\in 
\Omega_\nu\,\cap\Gamma.\tag{10}
\]
Now the domain $\Omega_\nu$ is bordered by the two simple curves
\[
\gamma_*=\{(x,y)\,\,\colon x=\rho_*(y)\}
\text{ and }
\gamma^*=\{(x,y)\,\,\colon x=\rho^*(y)\}
\]
where the inequalities $\alpha_*(\nu)\leq y\leq \alpha^*(\nu)$ hold.
Since $\partial_x H_k>\partial_x H_1$ we have
\[
H_1(\rho_*(y),y)<H_k(\rho^*(y),y),\quad 
\alpha_*(\nu)<y<\alpha^*(\nu).\tag{11}
\]
The subharmonic function
$V$ is equal to $H_1$ on $\gamma_*$ and it equals $H_k$ on
$\gamma^*$, so it follows from (11) that we must have
$V<H_k$ inside the domain
$\Omega_\nu\cap\Gamma$. This contradicts
equality (10) and the Sublemma is proved.\hfill $\Box$

\medskip

\noindent
\emph{Proof of Lemma 3.3, continued.}~The 
Sublemma shows that the locally piecewise harmonic function
$V$ inside $\Omega_\nu$
must
have at least one jump along $\Gamma\cap\Omega_\nu$, say from
from $H_j+c_j$ to $H_i+c_i$ for some indices $i<j$ in
$\{2,\ldots,k-1\}$. In other words, for some
$\alpha_*(\nu)<y_0<\alpha^*(\nu)$
there exists a small $\epsilon>0$ and constants $c_i,c_j$
such that
\[ 
V(\rho(y),y)=H_i(\rho(y),y)+c_i,\quad y'-\epsilon<y<y'\,,
\]
\[
V(\rho(y),y)=H_j(\rho(y),y)+c_j,\quad y'<y<y'+\epsilon\,.
\]

\noindent
By Proposition 2.3 and the strict monotonicity  of the sequence formed by
the $\partial_\Gamma$-derivatives of 
$H_2,\ldots,H_{k-1}$ in (ii) preceding the Non-return Lemma it follows that
\[ 
\partial_\Gamma H_j>\partial_\Gamma H_i
\]
on the whole of $\Gamma$. 
This is true for $\partial_\Gamma$-derivatives
whenever a jump occurs in some domain $\Omega_\nu$.
Hence, by the fact that the sequence in (ii) is strictly increasing
we cannot return to some $H_i$-function at a later stage if this function 
appears in some
$\Omega_\nu$-domain encountered previously.
Therefore $V|_{\Omega_\nu\cap\Gamma}$
can jump for at most $k-2$ values of $\nu$.
On the other hand, by the  Sublemma  a jump must always occur. We
conclude that the infinite sequence of intervals $\{J_\nu\}$
tending to $y=0$ as $\nu\to\infty$ cannot exist. 
This proves Lemma 3.3.\hfill $\Box$

\medskip

\noindent
{\bf 3.4.~Completing the proof of Proposition 3.1.}~Ignoring 
the case when $\rho_*(y)=\rho^*(y)$
in some interval $(0,\delta)$, in which case the equality 
$V=\text{max}(H_1,H_2)$
holds in a small rectangle 
\[
\square_0=\{(x,y)\,\colon\, -a<x<a,\, 0<y<b\},
\]
we have some positive $\delta_0$ from Lemma 3.3.
Set
\[
\Omega_0=\{(x,y)\,\colon\,  -a<x<\rho^*(y),\,0<y<\delta_0\}.
\]
In this domain $\partial V/\partial z$ is piecewise holomorphic 
with respect to
the $k-1$-tuple $g_1,\ldots,g_{k-1}$. By an induction over $k$
we may therefore assume that $V|_{\Omega_0}$
is locally piecewise harmonic with respect to $H_1,\ldots,H_{k-1}$ up to additive 
constants. 
Assume that $U_{k-1}\cap\Omega_0\neq\emptyset$.
Applying Theorem 1.2 it follows that
$U_{k-1}\cap\Omega_0=\mathfrak{s}^\uparrow(U_{k-1}\cap\Omega_0)$.
This gives a function $\rho_1(y)$ such that
\[ 
U_{k-1}\cap\Omega_0=\{(x,y)\,\colon\, x>\rho_1(y),\,0<y<\delta_0\}
\]
and there exists some $\delta_1\leq\delta_0$
such that
\[ 
\rho_*(y)\leq\rho_1(y)<\rho^*(y)\,,\quad\, 0<y<\delta_1.
\]
Since $U_{k-1}\cap\Omega_0$ is connected and $H_{k-1}$ only has to be 
determined up to a constant
we can assume that
$V=H_{k-1}$ in $U_{k-1}$ and then $H_{k-1}(\rho^*(y),y)=H_k(\rho^*(y),y)$
must hold when $0<y<\delta_1$, which  entails that
$H_{k-1}(0,0)=H_k(0,0)$.
For the next step we consider the domain
\[
\Omega_1=
\{(x,y)\,\colon\, -a<x<\rho_1(y),\,0<y<\delta_1\}\,.
\]
Again, if the closure of $U_{k-2}\cap\Omega_1$
contains the origin, then
it is equal to its forward star domain, which gives a function
$\rho_2(y)$ and some $0<\delta_2<\delta_1$
such that
 \[ 
 U_{k-2}\cap\Omega_1=
\{(x,y)\,\colon\, x>\rho_2(y),\,0<y<\delta_2\}\,.
\]
If it happens that the $g_{k-2}$-function is non-active, i.e., 
the closure of $U_{k-2}\cap \Omega_1$
is empty, we get a similar conclusion by taking the
largest 
integer $m\leq k-3$ such that the closure of $U_m\cap\Omega_1$
contains the origin.
We can continue in this way and arrive at the following result, 
where we use the notation
\[
\square^+(a,\delta)=\{(x,y)\,\colon\,  -a<x<a\,,\,0<y<\delta\}.
\]

\noindent {\bf 3.5.~Proposition}.~\emph{There exist  
a strictly increasing sequence
$1=j_1<\ldots <j_m=k$ and $a,\delta>0$ 
such that if for each $2\leq i\leq m-1$ we set
\[
\Omega_i=\bigl\{(x,y)\,\colon\, H_{j_i-1}(x,y)<
H_{j_i}(x,y)<H_{j_i+1}(x,y)\bigr\}\cap\square^+(a,\delta)
\]
then $V=H_{j_i}$ in $\Omega_i$
for $2<i\leq m-1$ while
$V=H_k$ when $H_k(x,y)>H_{j_{m-1}}(x,y)$
and $V=H_1(x,y)$ if $H_1(x,y)<H_{j_2}(x,y)$.}

\medskip

\noindent
{\bf Remark.}
A simpler way to express the above result 
is that in $\square^+(a,\delta)$ we have the
equality
$V=\text{max}(H_{j_1},\ldots,H_{j_m})$.
\medskip

\noindent
We can then proceed in exactly the same way in
the lower half-disk  where $y<0$
and obtain another $J$-sequence.
From this we conclude that if the resulting $\square$ 
(which is a neighborhood of the origin) is sufficiently small,
then
$U_\nu\cap\square$ has at most two connected components when
$2\leq \nu\leq k-1$ and is connected when $\nu=1$ or $k$.
This proves that $V$ is piecewise harmonic with respect to
$H_1,\ldots,H_k$ in $\square$, and 
completes the whole proof of Proposition 3.1.
\medskip

\noindent
{\bf 3.6.~Proof of Theorem 1.}~Given 
a domain $\Omega$
and some $k$-tuple $g_1,\ldots,g_k$ in
$\mathcal O(\Omega)$
we have the discrete set
\[
\sigma=\bigcap_{\nu\neq j}\,\{g_\nu=g_j\}\,.
\]

\noindent
If $z_0\in\Omega\setminus\sigma$
there exists some $\theta$ such that the sequence
$\{\mathfrak{Re}(e^{i\theta}g_\nu)(z_0)\}_1^k$
consists of distinct real numbers.
Up to a rotation we have the same local situation as in Proposition 3.1.
Hence $V$ is locally piecewise harmonic in
$\Omega\setminus\sigma$.
There remains only to study $V$ close to a single point
$z_0$ in
$\sigma$ and establish that it is  locally piecewise
harmonic in  a neighborhood of $z_0$.
Working locally we may take $z_0$ as the origin and in
a disk $D$ centered at $(0,0)$ we have the open subsets
$U_\nu=\{\partial V/\partial z=g_\nu\}$
of $D\setminus\text{supp}(\Delta(V))$.
Here the situation is more favorable than previously since we already know that
$V$ is locally piecewise harmonic in the punctured disk $\dot D$.
Moreover, to prove that $V$ is locally piecewise harmonic in
a neighborhood of the origin
it 
suffices to find some small $\delta>0$
such that the number of connected components of
each $U_\nu\cap D(\delta)$ is finite for all $1\leq\nu\leq k$, where 
$D(\delta)$ is the open disk of radius $\delta$ centered at $(0,0)$.
To achieve this  we will  decompose small discs into a finite
number of real analytic sectors
$\{\Omega_\alpha\}$ and  prove that
$U_\nu\cap\Omega_\alpha$ is empty or connected for each $\nu$ and $\alpha$.
For if this is  done then we may remove the union of
real analytic curves which border these sectors without
affecting the situation since this union is a null set, i.e., 
the locally integrable subharmonic function $V$ is not changed
by such a removal.
\medskip

\noindent
After these preliminary remarks we begin to 
construct suitable $\Omega$-sectors.
Consider  the harmonic functions
$\partial_x H_i-\partial_x H_\nu$ for pairs $i\neq\nu$.
Notice that such a function is identically zero if and only if
$H_i=H_\nu+cy$ for some constant $c$.
After a  rotation we may assume that this never happens and
get the real analytic set
\[
\Gamma=
\bigcup_{i\neq\nu}\,
\{\partial_x H_i=\partial_x H_\nu\}
\]
which  is described  in \S 2.5. So when $\delta>0$ is sufficiently small then
$D(\delta)\setminus\Gamma$
is a disjoint union of real analytic sectors
$\{\Omega_\alpha\}$.
It may occur that some sector contains a real line 
segment $0<x<\delta$ or $-\delta<x<0$.
Apart from this case a typical sector is given by
\[ \Omega=
\{(x,y)\colon\, \rho_*(y)<x<\rho^*(y),\, 0<y<\delta\}
\]
or by a similar sector in the lower half-disk where $-\delta<y<0$.
To handle  sectors that may potentially contain
a line segment on the $x$-axis we can simply
replace $x$ by
$y$  in the arguments above
and  start with the real analytic set
\[
\Gamma_1=
\bigcup_{i\neq\nu}\,
\{\partial_y(H_i)=\partial_y(H_\nu)\}\,.
\]
Then we again  obtain a finite number of real analytic sectors where
those which contain a line segment on the $x$-axis are 
defined by
\[
\{(x,y)\colon\, 0<x<\delta,\,\rho_*(x)<y<\rho^*(x)\}\,.
\]
Replacing $x$ by $y$ in Proposition 3.1 if necessary,
we conclude that  
the proof of Theorem 1 is finished if we can show
the following:

\medskip

\noindent
{\bf 3.7.~Proposition.}~\emph{Suppose 
$\partial_x H_1<\ldots<\partial_x H_k$ holds in $\Omega$, where 
$\Omega$ is a real analytic sector of the 
form}
\[
\{(x,y):\,\rho_*(y)<x<\rho^*(y),\, 0<y<\delta_0\}.
\] 
\emph{Then there exists $0<\delta<\delta_0$
such that if $\Omega(\delta)=\Omega\cap\{(x,y)\colon\, 0<y<\delta\}$ 
then $U_\nu\cap \Omega(\delta)$ is connected or empty
for every $\nu$.}

\medskip

\noindent
\emph{Proof.}~Arguments similar to those used in the proof of 
Proposition 3.1 yield the desired result.\hfill $\Box$

\medskip

\noindent
We have thereby completed the proof of Theorem 1.

\bigskip
\bigskip
\centerline{\bf\large
4.~On Algebraic Root Functions: Proof of Theorem 2}

\bigskip
\noindent
{\bf 4.1.~A  result inside sectors.}~As preparation for the proof of Theorem 2 we
first prove another result (Theorem 4.3 below) where
the harmonic functions under consideration 
are only defined in a real analytic sector.
Let
\[ 
\Omega=\{(x,y)\colon\,\rho_*(y)<x<\rho^*(y),\, 0<y<\delta_0\}
\]
be a real analytic sector and suppose that 
there are functions $\rho(y)<\rho_*(y)$ and
$\rho_1(y)>\rho^*(y)$ defining a larger sector
\[
\Omega^*=\{(x,y)\colon\,\rho(y)<x<\rho_1(y),\, 0<y<\delta_0\}.
\]
In $\Omega^*$ one is given a subharmonic function $V$
such that $\partial V/\partial z$ is piecewise
holomorphic with respect to some $k$-tuple $g_1,\ldots,g_k$ in
$\mathcal O(\Omega^*)$.
Since $\Omega^*$ is simply connected we have also a corresponding $k$-tuple of
harmonic functions $H_1,\ldots,H_k$ in $\Omega^*$.
Next,  \emph{assume} that Theorem 1 holds
in this situation, i.e., that $V$ is locally piecewise 
harmonic inside $\Omega^*$ with respect to the above 
$H$-functions plus constants.
In addition to this assumption we impose the following:
\medskip

\noindent
{\bf 4.2.~Condition on $\partial_\Gamma$-derivatives.}~For 
each pair $i\neq\nu$ and each constant $c$ set
\[
\Gamma(i,\nu,c)=\{
H_i-H_\nu=c\}\cap\,\Omega\,.
\]
For every such real analytic curve  we require that there exists
$\delta>0$ and
some permutation of $\{1,\ldots,k\}$ such that the inequalities 
\[
\partial_{\Gamma(i,\nu,c)}H_{j(1)}\leq
\ldots\leq
\partial_{\Gamma(i,\nu,c)}H_{j(k)}
\]
hold in  
\[
\Gamma_{(i,\nu,c)}\cap\,\{(x,y)\,\colon\, 0<y<\delta\}.
\]
Moreover, we require that there  exist index permutations so that
these inequalities hold for the tangential $H$-derivatives
along the two real analytic curves
$\{x=\rho_*(y)\}$ and $\{x=\rho^*(y)\}$.

\medskip

\noindent
{\bf 4.3.~Theorem.}~\emph{Under  the aforementioned conditions 
there exist  $\delta>0$ and
an increasing integer sequence
$1\leq j_1<\ldots<j_m\leq k$ such that} 
$V=\text{max}(H_{j_1},\ldots,H_{j_m})$ 
in
$\Omega\cap\{(x,y)\colon\, 0<y<\delta\}$.

\medskip

\noindent \emph{Proof.}~Follows by repeated use of Theorem 1 and arguments 
similar to those used in 
the proof of Proposition 3.1.\hfill $\Box$


\bigskip

\noindent
{\bf 4.4.}~As further preparation for the proof of
Theorem 2
we need some
results about root functions which arise as follows. Let
\[ 
f(z,y)=q_k(z)y^k+\ldots+q_1(z)y+q_0(z)
\]
be a polynomial in $y$ with coefficients
$q_\nu\in\mathcal O(D)$, where $D$ is an open disk centered at the origin.
We assume that $f$ has no multiple factors and get
the factorization
\[
f(z,y)=q_k(z)\cdot \prod_{\nu=1}^{k}\,
(y-\alpha_\nu(z)),
\]
where the $\alpha$-functions in general are multi-valued in the punctured disk
$\dot D$.
Set
\[
\Gamma=
\bigcup_{\nu\neq i}\,
\{\mathfrak{Re}(\alpha_i)-
\mathfrak{Re}(\alpha_\nu)=0\}\,.
\]
The real analytic set $\Gamma$ is to begin with only
defined in $\dot D$. Nevertheless, it extends to the whole disk
$D$ and becomes a union of smooth real analytic curves
passing through the origin. To see this we recall the classical
Normalisation Theorem saying that there exists an integer $M$ such that
if $\rho\colon\, \zeta\mapsto \zeta^M$ then 
$\alpha^*_\nu:=\alpha\circ\rho$
becomes meromorphic
in a disk of the $\zeta$-plane.
In this $\zeta$-disk we get the set
\[
\Gamma^*=
\bigcup_{\nu\neq i}\,
\{\mathfrak{Re}(\alpha^*_i)-
\mathfrak{Re}(\alpha^*_\nu)=0\}
\]
which is a disjoint union of smooth real analytic curves,  hence
so is the image $\Gamma=\rho(\Gamma^*)$.
Next we consider  the upper half-disk $D^+$
where $y>0$. Here we find single-valued branches of the root
functions and  consider their primitives
\[ 
A_\nu(z)=\int_p^z\,\alpha_\nu(w)dw,
\]
where the complex line integrals start from some
$p=ai$ with a small $a>0$.
In the Puiseux expansions of root functions
it may occur that $z^{-1}$ appears. So in $D^+$
we have
\[ 
A_\nu(z)= \lambda_\nu\cdot\text{Log}(z)+
\sum\, \psi_{i,\nu}(z)\cdot z^{\nu/M},
\] 
where the $\psi$-functions are meromorphic in $D$ and the $\lambda_\nu$'s 
are complex
numbers.
Note that any difference $A_i-A_\nu$ has a similar expansion.
Given some constant $c$
we use polar coordinates $z=re^{i\theta}$ to express a level curve as
\[
\Gamma=\{\mathfrak{Re}(A_i)-
\mathfrak{Re}(A_\nu)=c\}=
\{(r,\theta)\colon\,\,
u\cdot\text{Log}(r)-v\cdot\theta+\mathfrak{Re}(\Phi)(r,\theta)=c\},
\]
where $u,v$ are real constants and 
\[
\Phi(r,\theta)=
\sum_{j=0}^{M-1}\, r^{j/N}\cdot e^{j\theta/M}\cdot
\phi_\nu(r,\theta)
\]
with $\phi_0,\ldots,\phi_{M-1}$ meromorphic in $D$.
\medskip

\noindent
{\bf 4.5.~Tangential derivatives.}~In $D^+$ we get harmonic functions
$H_\nu=\mathfrak{Re}(A_\nu)$ which for each $\nu$ give
$\partial H_\nu/\partial z=\alpha_\nu$.
Along a level curve $\Gamma$ as above we consider a difference
\[
 \partial_\Gamma H_m
-\partial_\Gamma H_\ell\,,\quad
1\leq m,\ell\leq k\,.
\]
\noindent
Now we want to prove:

\medskip

\noindent {\bf 4.6.~Proposition.}~\emph{Unless
$\partial_\Gamma H_m
-\partial_\Gamma H_\ell$ is identically zero, there exists
$\delta>0$ such that this difference is
non-vanishing in $\Gamma(\delta):=\Gamma\cap D(\delta)$.}

\medskip

\noindent\emph{Proof.}~If $p\in\Gamma$ we notice that this
difference is zero at $p$
if and only if
\[ 
\mathfrak{Im}\!\left(\frac{\alpha_m-\alpha_\ell}{\alpha_i-\alpha_\nu}\right)
=0,\quad i\neq \nu.
\]

\noindent
The function $G:=\dfrac{\alpha_m-\alpha_\ell}{\alpha_i-\alpha_\nu}$
has a Puiseux series expansion in $D^+$:
\[
G(z)= \sum_{\nu=0}^{M-1}\,
g_\nu(z)\cdot z^{\nu/M}\,.
\]
Hence there only  remains to show:

\medskip

\noindent
\emph{Sublemma.}~\emph{There exists $\delta>0$ such that}
$p\in\Gamma(\delta)\implies\mathfrak{Im}(G)(p)\neq 0$.

\medskip

\noindent
{\em Proof.}~We use the existence of a holomorphic
map $\gamma$ from a complex $w$-disk onto $D$
such that
$ 
G\circ\gamma(w)=w^N
$
holds for some integer $N$.
Here 
$\{\mathfrak{Im}(G\circ\gamma)=0\}$
is a union of lines given by $\text{arg}(w)=m\pi/N\,,
\, 0\leq m\leq 2N-1$.
At the same time $\Gamma$
is the image of a curve
$\Gamma^*$ in the $w$-disk defined by an equation of the form 
\[
\Gamma^*=\{w \colon\,
\mathfrak{Re}[\lambda\cdot\text{Log}(w)+ S(w)]=c\},
\] 
where $S(w)$ is a meromorphic function.
In polar coordinates in the $w$-disk, $\Gamma^*$
is given by
\[ 
u\cdot\text{log}(r)+\mathfrak{Re}(S(re^{i\theta}))-v\theta-c=0,
\]
where $u,v$ are real constants.
The Sublemma follows since on each line in the zero set of
$G\circ\gamma$ where the $\theta$-angle is fixed, say $\theta=\theta_0$, 
it is obvious that
the function
\[ 
r\mapsto
u\cdot\text{log}(r)+\mathfrak{Re}(S(re^{i\theta_0}))-v\theta_0-c
\]
is non-vanishing for $0<r<\delta$
if $\delta$ is small enough, unless the function  happens to be identically 
zero.
This finishes the proof of Proposition 4.6.\hfill $\Box$

\medskip

\noindent
{\bf 4.7.~Proof of Theorem 2.~}Denote by
$\mathfrak{M}^+_{\text{alg}}$ the class of
probability measures $\mu$ such that
$\text{supp}(\mu)$ is a compact null set
and the Cauchy transform 
$\hat \mu(z)$ satisfies
an algebraic equation
\[
p_k(z)\cdot \hat \mu(z)^k+
\ldots+p_1(z)\cdot \hat \mu(z)+p_0(z)=0\,,\quad z\in 
{\bf{C}}\setminus
\text{supp}(\mu),
\]
where $p_0(z),\ldots,p_k(z)\in\mathbf{C}[z]$. 
Set $P(z,y)=p_k(z)\cdot y^k+
\ldots+p_1(z)\cdot y+p_0(z)$, which we assume to be irreducible
in ${\bf{C}}[z,y]$.
When the leading polynomial $p_k(z)\neq 0$ we have a factorization
\[
P(z,y)=p_k(z)\cdot
\prod_{i=1}^k\,
(y-\alpha_\nu(z))\,.
\]
We also get the rational  discriminant
\[
\mathfrak{D}(z)=\prod_{\nu\neq j}\,
(\alpha_\nu(z)-\alpha_j(z))\,.
\]
Let $\Sigma$ be the union of
$p_k^{-1}(0)$ and the zeros of
$\mathfrak{D}(z)$ in ${\bf{C}}\setminus p_k^{-1}(0)$.
Thus, if $U$ is a simply connected
subset of
${\bf{C}}\setminus\Sigma$
then the $k$-tuple of distinct
$\alpha$-roots are analytic functions in
$U$ and there exists some
$1\leq i\leq k$ such that
$\hat \mu=\alpha_i$ in $U$.

\medskip

\noindent
Consider now the subharmonic function
\[
V(z)=\int\,\text{log}\,|z-\zeta|\cdot d\mu(\zeta).
\]
Since $\partial V/\partial z=\hat \mu$ we can apply
Theorem 1 in the complement of
$\Sigma$. More  precisely, if $U$ as above is simply connected
we find harmonic functions
$H_1,\ldots,H_k$ in $U$ such that
$\partial_z(H_\nu)=\alpha_\nu$, $1\le \nu\le k$, 
and the restriction of $V|_U$ is locally
piecewise harmonic with respect to this $k$-tuple up to
additive constants.

\medskip

\noindent Next, using Proposition 4.6 we see that the conditions on
$\partial_\Gamma$-derivatives in \S 4.2
are satisfied when we consider suitable 
simply connected sectors
around each  individual point in $\Sigma$.
It follows again that the restriction of $V$ to each
such sector is
piecewise harmonic with respect to 
$H_1,\ldots,H_k$ up to constants. Applying Theorem 2.8 one 
finally arrives at Theorem 2.

\bigskip

\bigskip

\centerline{\bf\large
5.~Further Directions and Open Problems}
\bigskip

\noindent
{\bf 5.1.~Existence of measures in 
$\mathfrak{M}^+_{\text{alg}}$}.~Recall the class of probability measures 
$\mathfrak{M}^+_{\text{alg}}$ 
defined in \S 4.7. Consider a polynomial of the form
\[
P(y)= y+c_2y^2+\ldots+c_ky^k,\tag{12}
\]
where $k\geq 2$.
With $z$ as a new independent complex variable
we study the algebraic equation
\[ 
P(y)=\frac{1}{z}.
\]
\noindent
From (12) we see that if
$R$ is sufficiently large
then there exists a single-valued analytic function
$\alpha^*(z)$ defined in the exterior domain
$|z|>R$ whose Laurent expansion is
\[
\alpha^*(z)=\frac{1}{z}+a_2\cdot\frac{1}{z^2}+\ldots.
\]
\noindent
Let us assume that the zeros of $P(y)$ are simple.
In the complex $z$-plane we get the finite set
\[ 
\sigma=\left\{z=\frac{1}{P(\alpha)}\,\,\colon
P'(\alpha)=0\right\}.
\]
Clearly, $\sigma$ 
consists of $k-1$ points outside the orgin.
Now  $\alpha^*(z)$ extends to an (in general multi-valued) analytic function
defined in ${\bf{C}}\setminus \big(\sigma\cup\{0\}\big)$.
By an \emph{analytic tree} in
${\bf{C}}$ we mean a connected
set $\Gamma$ which is a finite union of
simple and closed real analytic Jordan arcs
and the open complement
${\bf{C}}\setminus\Gamma$
is connected. So by adding the point at infinity
the domain
$\Omega_\Gamma= \widehat{\bf{C}}\setminus\Gamma$ is 
simply connected. For
every such tree
$\Gamma$ which contains the set
 $\sigma\cup\{0\}$ the function $\alpha^*(z)$
extends from the exterior disk
$|z|>R$ to
a single-valued analytic function in
$\Omega_\Gamma$.
We also get the Riesz measure
$\mu_\Gamma$ supported by $\Gamma$ such that
\[
\alpha^*(z)=\hat\mu_\Gamma(z),
\quad z\in\Omega\setminus\Gamma.
\]
Since $\dfrac{1}{z}$ is the leading term in the 
Laurent expansion of $\alpha^*(z)$  we see that
\[
\int_\Gamma\, d\mu_\Gamma=1\,.\tag{13}
\]
The measure $\mu_\Gamma$ is in general complex-valued.
In fact, consider some relatively open Jordan arc
$\gamma\subset\Gamma$
which stays outside
$\sigma\cup\{0\}$.
Along the two opposite sides of
$\gamma$ we have two branches $\alpha_1(z)$ and $\alpha_2(z)$ of
$\alpha^*(z)$. By a classic formula  from analytic function theory
the restriction of $\mu_\Gamma$ to
the Jordan arc $\gamma$ is expressed by
\[
\frac{i}{\pi}\cdot [\alpha_2(z)-\alpha_1(z)]\cdot dz.
\]
To be precise, if $f(z)$ is a continuous function
whose compact support is disjoint from
$\Gamma\setminus\gamma$ then
\[ 
\int\, f\cdot d\mu=
i\cdot \int_\gamma\, f(z)\cdot [\alpha_2(z)-\alpha_1(z)]\cdot dz.
\]
\noindent
Notice that we can choose many different analytic 
trees
$\Gamma$ as above. For every such tree
the total variation of $\mu_\Gamma$ is $\geq 1$ by (13).
We propose the following:

\medskip

\noindent
{\bf Conjecture.}
\emph{There exists a unique analytic tree
$\Gamma$ such that
$\mu_\Gamma$ is a probability measure, i.e., 
$\mu_\Gamma\in \mathfrak{M}^+_{\text{alg}}$.}

\medskip

\noindent
{\bf Example.}~Consider the case $P(y)=y^2+y$.
Here $P'(y)= 2y+1$
and 
$\sigma=\{-4\}$.
Let $\Gamma$ be the (analytic) tree given by the real interval
$-4\leq x\leq 0$.
On this interval we define  the non-negative measure
\[ 
d\mu(x)=\frac{1}{2\pi}\cdot \frac{\sqrt{4+x}}{\sqrt{-x}}\,.
\]
Then we have
\[ 
\alpha^*(z)=
\hat\mu_\Gamma(z)
\]
in the complement of $\Gamma$.




\noindent
{\bf 5.2.~Combinatorics of subharmonic configurations.}~Given a $k$-tuple 
of harmonic functions $H_1,\ldots,H_k$ defined in an open connected set 
$\Omega\subset \mathbf{C}$ it is clear that there are locally only finitely 
many subharmonic configurations with respect to this
$k$-tuple.
It is natural to ask for the exact number of such configurations locally at a 
point $z\in \Omega$ in terms of the geometry of the convex hull of the 
gradients 
$\nabla H_\nu(z)$, $1\le \nu\le k$. Note that when all these 
gradients are extreme points of their convex hull Theorem 1.7 shows that 
there is only one possible configuration, namely the maximum of these 
harmonic functions. However, as seen in Example 2.11, 
in general there might be several such configurations. (Compare also 
Theorem 2.10.)


\medskip

\noindent
{\bf 5.3.~Plurisubharmonic configurations.}~An obvious question in this 
context is to try to extend some of our results to several variables, i.e., 
to study plurisubharmonic configurations with respect to given pluriharmonic 
functions. 

\medskip

\noindent
{\bf 5.4.~Configurations induced by fundamental solutions.}~A further 
interesting direction is to consider any partial differential operator of 
elliptic type for which one can define appropriate analogs of subharmonic 
functions. There are well known 
examples of possible such operators in the 
literature -- see, e.g., the subsolutions to
elliptic differential equations in \cite[Chapter 11]{J} 
as well as \cite{K,L}. It is 
known that the maximum of  subsolutions is a subsolution (cf., e.g., 
\cite{BS}), and hence one might for instance ask for conditions 
-- similar to the ones in Theorem 1.5 -- under
which this is the only subsolution.

\bigskip

\bigskip

\centerline{\bf\large Acknowledgements}

\bigskip

\noindent
The authors would like to thank Anders Melin at Lund University for insightful 
comments and inspirational discussions at an early stage of this work. The third author is also thankful to Boris Shapiro for discussions and suggestions. 


\end{document}